\documentclass[final,3p,times,numref]{elsarticle}

\usepackage{amsthm,amssymb,amsmath}
\usepackage{hyperref}

\newtheorem{theorem}{Theorem}[section]
\newtheorem{lemma}[theorem]{Lemma}

\newtheorem{conjecture}[theorem]{Conjecture}

\begin{document}

\begin{frontmatter}



\title{The Density and Quantitative Upper Bound of Integers Not Expressible as the Sum of a Prime and a Zumkeller Number} %% Article title


\author{Ting Hon Stanford Li} %% Author name

\affiliation{organization={Independent Mathematics Researcher},%Department and Organization
            country={Hong Kong}}

%% Abstract
\begin{abstract}
A Zumkeller number is a positive integer whose set of positive divisors can be partitioned into two disjoint subsets of equal sums. In this paper, we investigate the additive properties of these numbers in conjunction with prime numbers. Specifically, we study the natural density of the exceptional set \(\mathcal{E}\) consisting of all positive integers not expressible as the sum of a prime and a Zumkeller number. Let $\alpha(x) = \#\{n \le x : n \in \mathcal{E}\}$ denote its counting function. We establish a quantitative upper bound for $\alpha(x)$, demonstrating that $\alpha(x) \ll x(\log \log x)^{-A}$, where $A = \log_2\left(\frac{2\pi^2}{2\pi^2-3}\right) \approx 0.2378$, thereby proving that the natural density of this exceptional set is $0$.
\end{abstract}

%% Keywords
\begin{keyword}
Zumkeller numbers \sep Additive number theory \sep Natural density \sep Primes


\MSC[2020] 11B05 \sep 11P32 \sep 11B13

\end{keyword}

\end{frontmatter}



\section{Introduction}

    A Zumkeller number is a positive integer whose set of positive divisors can be partitioned into two disjoint subsets of equal sums. The concept was first introduced by LeVan in 1987, who named this class of integers integer-perfect numbers \cite{LeVan}. Subsequently, this class of integers was formalized and popularized by Zumkeller in 2003 via the On-Line Encyclopedia of Integer Sequences (sequence number A083207) \cite{OEIS}. Peng and Bhaskara Rao performed a systematic investigation of Zumkeller numbers, discovering many important properties \cite{Peng}.
    
    In a related direction, expressing integers as the sum of a prime and an element from a special arithmetic sequence is a classical theme in additive number theory. Romanoff proved that the set of integers expressible as the sum of a prime and a power of $2$ has a strictly positive lower asymptotic density \cite{Romanoff}, while Erd{\H{o}}s later proved the existence of an infinite arithmetic progression of odd integers that cannot be written in this form \cite{Erdos}. Recently, researchers have explored representations of integers as the sum of a prime and an element from other sequences, including Fibonacci numbers \cite{Lee, WangChen}, higher-order power sequences \cite{Brudern, Wang}, and general linear recurrences \cite{Dubickas, BallotLuca}.
    
    More recently, Somu, Kukla, and Tran \cite{Somu} studied the additive properties of Zumkeller numbers. They successfully proved that all sufficiently large integers can be written as the sum of two Zumkeller numbers, and similarly, as the sum of a Zumkeller number and a practical number. 
    
    Regarding prime-Zumkeller sums, they proved that there exist infinitely many positive integers that cannot be expressed as the sum of a prime and a Zumkeller number. Although this exceptional set is infinite, Somu, Kukla, and Tran conjectured that it has natural density zero.
    
    To formalize this conjecture, let $\mathcal{E}$ be the set of positive integers that are not expressible as the sum of a Zumkeller number and a prime:
    \[
    \mathcal{E} = \{ n \in \mathbb{Z}^+ : n \neq z + p, \text{ for any Zumkeller number } z \text{ and prime } p\}.
    \]
    For a real parameter $x \ge 1$, we define the corresponding counting function $\alpha(x)$ by 
    \[
    \alpha(x) = \# \{ n \le x : n \in \mathcal{E} \}.
    \]
    
    The natural density of a set $\mathcal{A} \subseteq \mathbb{Z}^+$, which we denote by $\mathbf{d}(\mathcal{A})$, is defined as the limit of the ratio of its counting function to $x$ as $x \to \infty$, provided the limit exists.
    
    These notations are used throughout the article.
    
    Using these notations, the conjecture of Somu, Kukla, and Tran can be stated as follows.
    
    \begin{conjecture}\label{con}
        The set $\mathcal{E}$ has natural density zero, that is,
        \[
        \mathbf{d}(\mathcal{E}) = \lim_{x \to \infty} \frac{\alpha(x)}{x} = 0.
        \]
    \end{conjecture}
    
    In this paper, we resolve Conjecture~\ref{con} by establishing a quantitative upper bound for the counting function of the exceptional set $\mathcal{E}$. Specifically, we prove the following main result:
    
    \begin{theorem}
        The counting function $\alpha(x)$ of the exceptional set $\mathcal{E}$ satisfies the asymptotic upper bound
        \[ \alpha(x) \ll \frac{x}{(\log \log x)^A}, \]
        where the exponent is given by
        \[ A = \log_2 \left( \frac{2\pi^2}{2\pi^2 - 3} \right) \approx 0.2378. \]
        Consequently, the natural density $d(\mathcal{E}) = 0$, thereby proving the conjecture of Somu, Kukla, and Tran.
    \end{theorem}
    
    The remainder of this paper is organized as follows. In Section 2, we present some important tools in analytic number theory and prior results on Zumkeller numbers. In Section 3, we establish three auxiliary lemmas related to prime-Zumkeller sums, and use them to prove our main theorem.

\section{Preliminaries}

    In this section, we recall several known results concerning Zumkeller numbers and classical analytic number theory. These will be used throughout the proof of our main theorem.
    
    We begin by recalling three established lemmas regarding Zumkeller numbers.
    
    The first two lemmas were established by Peng and Bhaskara Rao.
    
    \begin{lemma}[{\cite[Proposition~9]{Peng}}]\label{zumkeller_practical}
        For any prime $p \neq 2$ and positive integer $k$ with $p \leq 2^{k+1} - 1$, the number $2^kp$ is a Zumkeller number.
    \end{lemma}
    
    \begin{lemma}[{\cite[Corollary~5]{Peng}}]\label{zumkeller}
        If the positive integer $n$ is Zumkeller and $w$ is a positive integer which is relatively prime to $n$, then $nw$ is Zumkeller.
    \end{lemma}
    
    The third lemma was proven by Somu, Kukla, and Tran.
    
    \begin{lemma}[Adapted from Somu, Kukla, and Tran \cite{Somu}, Lemma 3.3]\label{zumkeller2}
        Let $p_i$ denote the $i$-th prime, then
        \begin{enumerate}
            \item For each positive integer $k$, there exists a unique smallest index $f(k)$ such that
            \[
            \sigma(p_k p_{k+1} \cdots p_{f(k)}) = (1 + p_k)(1 + p_{k+1}) \cdots (1 + p_{f(k)}) \geq 2p_k p_{k+1} \cdots p_{f(k)}.
            \]
            \item Let $A_k = p_k p_{k+1} \cdots p_{f(k)}$. Then there exists a positive integer threshold $k_1 \geq 2$ such that $A_k$ is an odd Zumkeller number for all $k > k_1$.
        \end{enumerate}
    \end{lemma}
    
    We also require two classical analytic number theory results, namely Linnik's theorem and Mertens' Second Theorem, for proving our main theorem.
    \begin{lemma}[Linnik's theorem, \cite{Linnik}]\label{Linnik}
        Let $a$ and $d$ be positive, coprime integers such that $1 \le a < d$. 
        The smallest prime $p$ satisfying $p \equiv a \pmod d$ also satisfies 
        \[
        p \ll d^L,
        \]
        where $L$ is the Linnik's constant.
    \end{lemma}
    
    \begin{lemma}[Mertens' Second Theorem, {\cite[Theorem~4.12]{Apostol}}]\label{merten}
        Let $x \ge 2$. The sum of the reciprocals of all prime numbers $p$ less than or equal to $x$ satisfies the asymptotic formula:
        \[
        \sum_{p \le x} \frac{1}{p} = \log(\log x) + M + O\left(\frac{1}{\log x}\right)
        \]
        where $M \approx 0.261497$ is the Meissel-Mertens constant.
    \end{lemma}

\section{Proof of the Main Result}

    In this section, we establish three auxiliary lemmas that serve as the foundation for our main theorem.
    
    \begin{lemma}\label{odd}
        Let $p$ be an odd prime. Let $t_p$ be the smallest integer such that $2^{t_p+1} -1 > p$ and $k_p = p2^{t_p}$. Then there exists a positive integer $N_p \ll p^{3L}$ such that all odd integers $n > N_p$ satisfying $n \not\equiv 0 \pmod p$ can be expressed as the sum of a prime and a Zumkeller number which is a multiple of $k_p$.
    \end{lemma}
    
    \begin{proof}
        By Lemma~\ref{zumkeller_practical}, $k_p$ is a Zumkeller number. Also, note that for all positive integers $m$, $k_p$, which is a multiple of $p$, is always coprime with $2pm+1$. So by Lemma~\ref{zumkeller}, $k_p(2pm+1)$ is a Zumkeller number.
        
        Our aim is to express an arbitrary odd positive integer $n$ in the following form:
        \[ n = g_n + k_p(2pm+1), \]
        where $g_n$ is a prime and $m$ is a positive integer.
        
        Define $y_p = 2pk_p = 2p^2 2^{t_p}$. Rearranging our target equation, we require a prime $g_n$ that satisfies:
        \[ n - g_n = 2pk_pm + k_p = y_pm + k_p, \]
        which satisfies the congruence:
        \[ g_n \equiv n - k_p \pmod{y_p}. \]
        
        By Dirichlet's Theorem on Arithmetic Progressions, a prime $g_n$ exists in this residue class if and only if $\gcd(n - k_p, y_p) = 1$. The only prime factors of $y_p = 2p^2 2^{t_p}$ are $2$ and $p$. We consider the coprimality condition across both factors:
        \begin{enumerate}
            \item Modulo 2: The integer $n$ is odd, and $k_p = p2^{t_p}$ is even since $t_p \ge 1$. The difference $n - k_p$ is therefore always odd, so $\gcd(n - k_p, 2) = 1$.
            \item Modulo p: Since $k_p \equiv 0 \pmod p$, we have $n - k_p \equiv n \pmod p$. By hypothesis, $n \not\equiv 0 \pmod p$, and so $\gcd(n - k_p, p) = 1$.
        \end{enumerate}
        
        Thus, for every odd integer $n \not\equiv 0 \pmod p$, the residue class $r_n \equiv n - k_p \pmod{y_p}$ is coprime to $y_p$. 
        
        For each coprime residue class $r \pmod{y_p}$, let $q(r)$ denote the smallest prime congruent to $r$ modulo $y_p$. The existence of these primes is guaranteed by Dirichlet's Theorem on Arithmetic Progressions. Since the size of the set of coprime residue classes modulo $y_p$ is finite, we can define a threshold $N_p$ as follow:
        \[
        N_p = \max_{\gcd(r, y_p)=1} \{q(r)\} + k_p.
        \]
        
        We first show that $N_p \ll p^{3L}$. By Linnik's theorem, the smallest prime $q(r) \equiv r \pmod{y_p}$ satisfies $q(r) \ll y_p^{\,L}$. 
        By the definition of $t_p$, we have $2^{t_p} \asymp p$, and therefore,
        \[
        y_p = 2p^22^{t_p} \ll p^3.
        \]
        
        Hence,
        \[
        q(r) \ll y_p^{\,L} \ll (p^3)^L = p^{3L}.
        \]

        Since $k_p = p2^{t_p} \ll p^2  \ll p^{3L}$ for any Linnik constant $L>1$, adding $k_p$ does not affect the exponent. We obtain the following quantitative upper bound for $N_p$:
        \[
        N_p = \max_{\gcd(r, y_p)=1} \{q(r)\} + k_p \ll p^{3L}.
        \]
        
        Now, let $n > N_p$ be any odd integer satisfying $n \not\equiv 0 \pmod p$, and let $r_n \equiv n - k_p \pmod{y_p}$ be its corresponding coprime residue class. We choose $g_n = q(r_n)$. Because $n > N_p$, it follows that 
        \[
        n - k_p > \max_{\gcd(r, y_p)=1} \{q(r)\} \geq q(r_n) = g_n,
        \]
        which guarantees that $n - g_n - k_p > 0$. Furthermore, since $g_n \equiv r_n \equiv n - k_p \pmod{y_p}$, the difference $n - g_n - k_p$ is a strictly positive multiple of $y_p$. Consequently, 
        \[
        m = \frac{n - g_n - k_p}{y_p}
        \]
        is a strictly positive integer. Rearranging this expression yields our target representation:
        \[
        n = g_n + k_p(2pm + 1).
        \]
        Since $g_n$ is prime and $k_p(2pm + 1)$ is a Zumkeller number by Lemma~\ref{zumkeller}, it follows that $n$ can be expressed as the sum of a prime and a Zumkeller number, which completes the proof.
    \end{proof}
    
    \begin{lemma}\label{even}
        Let $z$ be an odd Zumkeller number. Let $v_{1}, v_{2}, \dots, v_{i}$ be all the distinct prime divisors of $z$. Then there exists a positive integer $N' \ll z^{2L}$ such that all even integers $n > N'$ satisfying $n \not\equiv 0 \pmod{v_{j}}$ for all $j \in \{1, \dots, i\}$ can be expressed as the sum of a prime and a Zumkeller number which is a multiple of $z$.
    \end{lemma}
    
    \begin{proof}
        Let $P = \prod_{j=1}^i v_{j}$ denote the product of all distinct prime factors of $z$. 
        
        Since for all positive integer $m$, $2Pm + 1 \equiv 1 \pmod{v_{j}}$ for every prime divisor $v_{j}$ of $z$, it follows that $\gcd(2Pm + 1, z) = 1$. Because $z$ is a Zumkeller number and $2Pm + 1$ is coprime to it, by Lemma~\ref{zumkeller}, the product $z(2Pm + 1)$ is a Zumkeller number which is a multiple of $z$.
        
        We aim to find an additive representation of an arbitrary even integer $n$ of the form:
        \[ n = g_n + z(2Pm + 1), \]
        where $g_n$ is a prime and $m$ is a positive integer. We define $y' = 2zP$. Rearranging this target equation yields:
        \[ n - g_n = 2zPm + z = y'm + z. \]
        Thus, the prime $g_n$ must satisfy the congruence:
        \[ g_n \equiv n - z \pmod{y'}. \]
        
        By Dirichlet's Theorem on Arithmetic Progressions, a prime $g_n$ is guaranteed to exist in this residue class if and only if $\gcd(n - z, y') = 1$. The set of prime factors of $y' = 2zP$ consists of $2$ and the distinct prime factors of $z$, namely $\{2, v_{1}, \dots, v_{i}\}$. Evaluating the coprimality condition modulo each prime factor, we find:
        \begin{enumerate}
            \item Modulo 2: The integer $n$ is even, and $z$ is odd. The difference $n - z$ is therefore always odd, so $\gcd(n - z, 2) = 1$.
            \item Modulo $v_j$: Since $z \equiv 0 \pmod{v_j}$, we have $n - z \equiv n - 0 \equiv n \pmod{v_{j}}$. Therefore, to ensure $\gcd(n - z, v_j) = 1$, the integer $n$ must satisfy $n \not\equiv 0 \pmod{v_{j}}$ for all $j \in \{1, \dots, i\}$, which is guaranteed by our hypothesis.
        \end{enumerate}
        
        For each coprime residue class $r \pmod{y'}$, let $q(r)$ denote the smallest prime congruent to $r$ modulo $y'$. Since the set of coprime residue classes modulo $y'$ is finite, we can define a threshold $N'$ by taking the maximum over this finite set:
        \[
        N' = \max_{\gcd(r, y')=1} \{q(r)\} + z.
        \]
        
        We first show that $N' \ll z^{2L}$. Since $P$ is the product of the distinct prime divisors of $z$, we have $P \le z$, and therefore,
        \[
        y' = 2zP \le 2z^2.
        \]

        By Linnik's theorem, the smallest prime $q(r) \equiv r \pmod{y'}$ satisfies $q(r) \ll (y')^{L} \ll (z^2)^L = z^{2L}$. Because $z \ll z^{2L}$ for any Linnik constant $L>1$, adding the term $z$ does not affect the exponent. Since the number of coprime residue classes modulo $y'$ is finite, the maximum defining $N'$ is well-defined, and we obtain:
        \[
        N' = \max_{\gcd(r, y')=1} \{q(r)\} + z \ll z^{2L}.
        \]
        
        Now, let $n > N'$ be any even integer satisfying $n \not\equiv 0 \pmod{v_j}$ for all $j \in \{1, \dots, i\}$, and let $r_n \equiv n - z \pmod{y'}$ be its corresponding coprime residue class. We choose $g_n = q(r_n)$. Because $n > N'$, it follows that 
        \[
        n - z > \max_{\gcd(r, y')=1} \{q(r)\} \geq q(r_n) = g_n,
        \]
        which guarantees that $n - g_n - z > 0$. Furthermore, since $g_n \equiv r_n \equiv n - z \pmod{y'}$, the difference $n - g_n - z$ is a strictly positive multiple of $y'$. Consequently, 
        \[
        m = \frac{n - g_n - z}{y'}
        \]
        is a strictly positive integer. Rearranging this expression yields our target representation:
        \[
        n = g_n + z(2Pm + 1).
        \]
        Since $g_n$ is prime and $z(2Pm + 1)$ is a Zumkeller number by Lemma~\ref{zumkeller}, it follows that $n$ can be expressed as the sum of a prime and a Zumkeller number that is a multiple of $z$, which completes the proof.
    \end{proof}
    
    \begin{lemma}\label{zumkeller3}
        For all positive integers $r$, there exist odd Zumkeller numbers $z_1, z_2, \dots, z_r$ which are pairwise coprime and $\frac{\phi(z_t)}{z_t} > \frac{3}{2\pi^2}$ for $t=1,2,\dots,r$. Furthermore, $z_r \ll \exp\left(\exp(D \cdot 2^r)\right)$ for a positive constant $D$.
    \end{lemma}
    
    \begin{proof}
        By Lemma~\ref{zumkeller2}, there exists an integer $k_1$ such that for any starting prime index $k > k_1$, the product of consecutive primes $A_k = p_k p_{k+1} \cdots p_{f(k)}$ is an odd Zumkeller number. Here, $f(k)$ is defined as the smallest index where $\sigma(p_k \cdots p_{f(k)}) = (1 + p_k)\cdots(1 + p_{f(k)}) \ge 2p_k \cdots p_{f(k)}$, its existence is guaranteed by Lemma~\ref{zumkeller2}.
        
        By the minimality definition of $f(k)$, we have:
        \[ \prod_{j=k}^{f(k)-1}\left(1 + \frac{1}{p_j}\right) < 2 \implies \prod_{j=k}^{f(k)}\left(1 + \frac{1}{p_j}\right) < 2\left(1 + \frac{1}{p_{f(k)}}\right)<2(1+1)=4.\]
        
        Moreover, by the Euler Product Formula (see, for example, \cite[Theorem 11.7]{Apostol}), we have:
        \[ \prod_{i=1}^{\infty} \left(1 - \frac{1}{p_i^2}\right) = \frac{6}{\pi^2}  \implies \prod_{j=k}^{f(k)}\left(1 - \frac{1}{p_j^2}\right) > \frac{6}{\pi^2}.\]
        
        Combining, we have:
        \[
        \frac{\phi(A_k)}{A_k}=\prod_{j=k}^{f(k)}\left(1 - \frac{1}{p_j}\right)=\frac{\prod_{j=k}^{f(k)}\left(1 - \frac{1}{p_j^2}\right)}{\prod_{j=k}^{f(k)}\left(1 + \frac{1}{p_j}\right)} >\frac{\frac{6}{\pi^2}}{4} = \frac{3}{2\pi^2}.
        \]
        
        We define $s_1=k_1+1$; and for all positive integers $t>1$, we define $s_t=f(s_{t-1})+1$. Let $z_t=A_{s_t}=p_{s_t}p_{s_t+1}p_{s_t+2}\dots p_{f(s_t)}$ for $t=1,2,\dots,r$. Then $z_1, z_2, \dots, z_r$ are odd Zumkeller numbers which are pairwise coprime. With the above argument, we have proved that $\frac{\phi(z_t)}{z_t} > \frac{3}{2\pi^2}$ for $t=1,2,\dots,r$.
        
        We then prove that $z_r \le \exp\left(\exp(D \cdot 2^r)\right)$ for a positive constant $D$.
        From the minimality definition of $f(k)$, we have established that:
        \[
        \prod_{j=k}^{f(k)}\left(1 + \frac{1}{p_j}\right) < 2\left(1 + \frac{1}{p_{f(k)}}\right) < 2\left(1 + \frac{1}{p_k}\right).
        \]
        Taking the natural logarithm of both sides and applying the inequality $\log(1+x)<x$ for all $x>0$, we have:
        \[
        \sum_{j=k}^{f(k)} \log\left(1 + \frac{1}{p_j}\right) < \log 2 + \log\left(1 + \frac{1}{p_k}\right) \implies \sum_{j=k}^{f(k)} \log\left(1 + \frac{1}{p_j}\right)<\log 2+\frac{1}{p_k}.
        \]
        Applying the first-order Taylor series expansion $\log(1 + x) = x + O(x^2)$ for small $x = \frac{1}{p_j}$, we obtain:
        \[
        \sum_{j=k}^{f(k)} \frac{1}{p_j} = \log 2 + O\left(\frac{1}{p_k}\right).
        \]
        By applying Mertens' Second Theorem (Lemma~\ref{merten}), we have:
        \[
        \left( \log\log p_{f(k)} + M + O\left(\frac{1}{\log p_{f(k)}}\right) \right) - \left( \log\log p_k + M + O\left(\frac{1}{\log p_k}\right) \right) = \log 2 + O\left(\frac{1}{p_k}\right),
        \]
        where $M$ is the Meissel-Mertens constant. Simplifying the equation yields:
        \[
        \log\left(\frac{\log p_{f(k)}}{\log p_k}\right) = \log 2 + O\left(\frac{1}{\log p_k}\right).
        \]
        As a result,
        \[
        \frac{\log p_{f(k)}}{\log p_k} = 2 \left(1 + O\left(\frac{1}{\log p_k}\right)\right) \implies \log p_{f(k)} = 2\log p_k + O(1),
        \]
        This implies that $p_{f(k)} \asymp p_k^2$. Therefore, we obtain:
        \[
        \log p_{s_t} = \log p_{f(s_{t-1})+1} \sim \log p_{f(s_{t-1})} \approx 2\log p_{s_{t-1}}.
        \]
        Iterating this from $t = 1$ up to $t = r$ yields:
        \[
        \log p_{s_r} \ll 2^r \implies p_{s_r} \ll \exp(C \cdot 2^r),
        \]
        for a positive constant $C$.
        
        Now, we approximate the magnitude of $z_r$ by applying Chebyshev's $\theta$-function estimate $\theta(x) = \sum_{p \le y} \log p \sim y$ (see, e.g., \cite[Theorem 4.4]{Apostol}):
        \begin{align*}
        \log z_r &= \sum_{j=s_r}^{f(s_r)} \log p_j = \theta(p_{f(s_r)}) - \theta(p_{s_r}-1) \\
        &\sim p_{f(s_r)} - p_{s_r} \\
        &\asymp p_{s_r}^2 - p_{s_r} \\
        &\asymp p_{s_r}^2 \\
        &\ll \left[\exp(C \cdot 2^r)\right]^2 \\
        &= \exp(2C \cdot 2^r) \\
        &= \exp(D \cdot 2^r),
        \end{align*}
        for a positive constant $D$.
        
        Exponentiating the expression gives us:
        \[
        z_r \ll \exp\left(\exp(D \cdot 2^r)\right).
        \]
    \end{proof}
    
    We are ready to prove the main theorem.
    
    \begin{proof}[Proof of Theorem 1.2]
    
        We partition the set $\mathcal{E}$ into its odd and even parity components, letting $\mathcal{E}_{\text{odd}} = \mathcal{E} \cap (2\mathbb{Z}^+ - 1)$ and $\mathcal{E}_{\text{even}} = \mathcal{E} \cap 2\mathbb{Z}^+$. For any subset $\mathcal{A} \subseteq \mathbb{Z}^+$, we let $|\mathcal{A}(x)| = \#\{n \le x : n \in \mathcal{A}\}$ denote its corresponding counting function up to $x$. It follows that $\alpha(x) = |\mathcal{E}_{\text{odd}}(x)| + |\mathcal{E}_{\text{even}}(x)|$.
        
        \smallskip
        \noindent\textbf{Case 1: The Odd Integers ($\mathcal{E}_{\text{odd}}$).} 
        Suppose $x$ is sufficiently large. We choose an odd prime $p \asymp x^{\frac{1}{3L+1}}$. By Lemma~\ref{odd}, there exists a positive integer threshold $N_p \ll p^{3L}$ such that every odd integer $n > N_p$ satisfying $n \not\equiv 0 \pmod p$ is expressible as the sum of a prime and a Zumkeller number. 
        
        Consequently, any unexpressible odd integer $n \in \mathcal{E}_{\text{odd}}$ up to $x$ must either be smaller than $N_p$ or be a multiple of $p$. We count these elements:
        \[
        |\mathcal{E}_{\text{odd}}(x)| \le N_p + |\{n \le x : n \equiv 1 \pmod 2 \text{ and } n \equiv 0 \pmod p\}| \le N_p + \frac{x}{2p}.
        \]
        Since we choose $p \asymp x^{\frac{1}{3L+1}}$, we can substitute our estimate $N_p \ll p^{3L}$ and obtain:
        \[
        |\mathcal{E}_{\text{odd}}(x)| \ll p^{3L} + \frac{x}{p} \ll (x^{\frac{1}{3L+1}})^{3L} = x^{1-\frac{1}{3L+1}}.
        \]
        
        \smallskip
        \noindent\textbf{Case 2: The Even Integers ($\mathcal{E}_{\text{even}}$).} 
        
        Let $z_1, z_2, \dots, z_r$ be odd Zumkeller numbers which are pairwise coprime and $\frac{\phi(z_t)}{z_t} > \frac{3}{2\pi^2}$ for $t=1,2,\dots,r$. By Lemma~\ref{zumkeller3}, they exist. Let $M_t$ be the set of all Zumkeller numbers which are multiples of $z_t$. Also, let $N'_1, N'_2, \dots, N'_r$ be the respective positive integers defined in Lemma~\ref{even}.
        
        Let us evaluate the natural density of integers that cannot be expressed as the sum of a prime and an odd Zumkeller number in the set $M_1 \cup M_2 \cup M_3 \cup \dots \cup M_r$:
        
        By Lemma~\ref{even}, for each $t \in \{1, 2, \dots, r\}$, any even integer $n > N'_t$ is guaranteed to be expressible as a sum of a prime and an odd Zumkeller number in $M_t$ if it satisfies $\gcd(n, z_t) = 1$. Therefore, for any even integer $n > \max \{N'_1, N'_2, \dots, N'_r \}$, the condition for $n$ being not expressible by any of the sets $M_t$ requires that $n$ fails the coprimality condition for every single $z_t$, meaning $\gcd(n, z_t) \neq 1$ for all $t \in \{1, 2, \dots, r\}$.
        
        For each $t$, the density of integers satisfying $\gcd(n,z_t) \neq 1$ is $1-\frac{\phi(z_t)}{z_t}$. Because the Zumkeller numbers $z_1, z_2, \dots, z_r$ are pairwise coprime, their sets of prime factors are disjoint. Furthermore, since each $z_t$ is strictly odd, their prime divisors exclude the prime $2$. By the Chinese Remainder Theorem, the parity condition $n \equiv 0 \pmod 2$ is mutually independent of the congruence conditions modulo the odd prime factors of $z_t$.
        
        Recall that any unexpressible even integer $n \in \mathcal{E}_{\text{even}}$ up to $x$ must either be smaller than $\max\{N'_1,N'_2,\dots,N'_r\}$ or be a multiple of $z_1, z_2, \dots, z_r$. Moreover, by Lemma~\ref{even}, $\max\{N'_1,N'_2,\dots,N'_r\} \ll z_r^{2L}$. Overall, we obtain:
        \[
        |\mathcal{E}_{\text{even}}(x)| \ll z_r^{2L} + \frac{x}{2} \prod_{t=1}^r \left(1 - \frac{\phi(z_t)}{z_t}\right).
        \]
        Applying the bound $\frac{\phi(z_t)}{z_t} > \frac{3}{2\pi^2}$ across all $t$ simplifies the relationship as:
        \[
        |\mathcal{E}_{\text{even}}(x)| \ll z_r^{2L} + x \left(\frac{2\pi^2-3}{2\pi^2}\right)^r.
        \]
        
        Suppose $x$ is sufficiently large, we choose $r$ such that $z_r^{2L} \le x^{\frac{1}{2}}$. By Lemma~\ref{zumkeller3}, $z_r$ satisfies $\log z_r \ll \exp(D \cdot 2^r)$. Taking the natural logarithm of both sides of $z_r^{2L} \le x^{\frac{1}{2}}$ yields:
        \[
        2L \log z_r \le \frac{1}{2} \log x \implies \log z_r \ll \log x \implies \exp(D \cdot 2^r) \ll \log x \implies 2^r \ll \log \log x.
        \]
        
        We now apply this quantitative upper bound and obtain:
        \begin{align*}
        x \left(\frac{2\pi^2-3}{2\pi^2}\right)^r &= x \left(2^{\log_2\left(\frac{2\pi^2-3}{2\pi^2}\right)}\right)^r \\
        &= x \left(2^r\right)^{\log_2\left(\frac{2\pi^2-3}{2\pi^2}\right)} \\
        &\ll x (\log \log x)^{\log_2\left(\frac{2\pi^2-3}{2\pi^2}\right)} \\
        &= \frac{x}{(\log \log x)^{-\log_2\left(\frac{2\pi^2-3}{2\pi^2}\right)}} \\
        &= \frac{x}{(\log \log x)^A},
        \end{align*}
        where $A = \log_2\left(\frac{2\pi^2}{2\pi^2-3}\right) \approx 0.2378$.
        
        Since the term $z_r^{2L} \le x^{1/2}$ decays significantly faster than $\frac{x}{(\log \log x)^A}$, we have:
        \[
        |\mathcal{E}_{\text{even}}(x)| \ll z_r^{2L} + x \left(\frac{2\pi^2-3}{2\pi^2}\right)^r\ll \frac{x}{(\log \log x)^A}.
        \]
        
        \smallskip
        \noindent\textbf{Conclusion.} 
        Combining both cases, the even case dominates. Consequently,
        \[
        \alpha(x) = |\mathcal{E}_{\text{odd}}(x)| + |\mathcal{E}_{\text{even}}(x)| \ll x^{1-\frac{1}{3L+1}} + \frac{x}{(\log \log x)^A} \ll \frac{x}{(\log \log x)^A}.
        \]
        Dividing by $x$ and taking the limit as $x \to \infty$ yields $\lim_{x \to \infty} \frac{\alpha(x)}{x} = 0$, which establishes that $d(\mathcal{E}) = 0$.
    \end{proof}

%%       \citet{<label>} ==> Lamport (1994)
%%       \citep{<label>} ==> (Lamport, 1994)
%%


\begin{thebibliography}{0}

\bibitem[Apostol(1976)]{Apostol}
  T. M. Apostol,
  \textit{Introduction to Analytic Number Theory},
  Springer-Verlag, New York,
  1976.

\bibitem[Dubickas(2013)]{Dubickas}
  A. Dubickas,
  Sums of primes and quadratic linear recurrence sequences,
  \textit{Acta. Math. Sin.-English Ser.} \textbf{29} (2013) 2251--2260,
  \url{https://doi.org/10.1007/s10114-013-1059-x}.

\bibitem[Ballot and Luca(2013)]{BallotLuca}
  C. Ballot and F. Luca,
  On the sumset of the primes and a linear recurrence,
  \textit{Acta Arith.} \textbf{161} (2013) 33--46,
  \url{https://doi.org/10.4064/aa161-1-3}.

\bibitem[Br\"udern and Perelli(1996)]{Brudern}
  J. Br\"udern and A. Perelli,
  The addition of primes and power,
  \textit{Canad. J. Math.} \textbf{48}(3) (1996) 512--526,
  \url{https://doi.org/10.4153/CJM-1996-026-8}.

\bibitem[Erd\H{o}s(1950)]{Erdos}
  P. Erd\H{o}s,
  On integers of the form $2^k + p$ and some related problems,
  \textit{Summa Brasil. Math.} \textbf{2} (1950) 113--123.

\bibitem[Lee(2010)]{Lee}
  K. S. Enoch Lee,
  On the sum of a prime and a Fibonacci number,
  \textit{Int. J. Number Theory} \textbf{6} (7) (2010) 1669--1676,
  \url{https://doi.org/10.1142/S1793042110003708}.

\bibitem[LeVan(1987)]{LeVan}
  M. O. LeVan,
  Integer-perfect numbers,
  \textit{Journal of Natural Sciences and Mathematics} \textbf{27}(2) (1987) 33--50.

\bibitem[Linnik(1944)]{Linnik}
  U. V. Linnik,
  On the least prime in an arithmetic progression. I. The basic theorem,
  \textit{Rec. Math. [Mat. Sbornik] N.S.} \textbf{15}(57) (1944) 139--178.

\bibitem[OEIS Foundation Inc.(2026)]{OEIS}
  OEIS Foundation Inc.,
  The On-Line Encyclopedia of Integer Sequences,
  \url{https://oeis.org}.

\bibitem[Peng and Rao(2013)]{Peng}
  Y. Peng and K. P. S. Bhaskara Rao,
  On Zumkeller numbers,
  \textit{J. Number Theory} \textbf{133}(4) (2013) 1135--1155,
  \url{https://doi.org/10.1016/j.jnt.2012.09.020}.

\bibitem[Romanoff(1934)]{Romanoff}
  N. P. Romanoff,
  \textsl{\"{U}ber einige S\"{a}tze der additiven Zahlentheorie},
  \textit{Math. Ann.} \textbf{109} (1934) 668--678,
  \url{https://doi.org/10.1007/BF01449161}.

\bibitem[Somu et al.(2023)]{Somu}
  S. T. Somu, A. Kukla and D. V. K. Tran,
  Some results on Zumkeller numbers,
  preprint arXiv:2310.14149 [math.NT], 2023.

\bibitem[Wang and Chen(2023)]{WangChen}
  R.-J. Wang and Y.-G. Chen,
  On the sum of a Fibonacci number and a prime,
  \textit{Int. J. Number Theory} \textbf{19}(4) (2023) 873--889,
  \url{https://doi.org/10.1142/S1793042123500446}.

\bibitem[Wang(2012)]{Wang}
  Y. C. Wang,
  On the sum of a prime and a $k$-th power of prime in short intervals,
  \textit{Acta Math. Hungar.} \textbf{135} (2012) 248--269,
  \url{https://doi.org/10.1007/s10474-011-0193-0}.

\end{thebibliography}


\section{Declaration of generative AI and AI-assisted technologies in the manuscript preparation process}

During the preparation of this work, the author used Gemini 3 for assisting with manuscript formatting and language polishing. No mathematical content within this paper was generated by generative AI. All mathematical content is the original work of the author. After using this tool, the author reviewed and edited the content as needed and takes full responsibility for the content of the published article.

\end{document}